\documentclass[smallextended,referee,envcountsect,]{svjour3}
\smartqed
\usepackage{graphicx}

\usepackage{amssymb,amsmath}
\usepackage{hyperref}

\usepackage{xcolor}

\def\R{\mathbb{R}}
\def\la{\langle}
\def\ra{\rangle}
\def\pol{\operatorname{Polar}}

\def\resid{\epsilon_{\mathrm{cert}}}
\def\resar{\left(\xi \mid P_\tau, Q_1\right)}
\def\Opt{\mathrm{Opt}}
\def\dom{\mathrm{dom}}
\renewcommand\int{\operatorname{int}}

\def\nul{{\mathbf 0}}
\def\last{\left[\begin{smallmatrix} \nul \\ 1 \end{smallmatrix}\right]}

\def\It{\mathcal{I}_\tau}
\def\Pt{\mathcal{P}_\tau}
\def\Nt{\mathcal{N}_\tau}

\begin{document}

\title{Accuracy Certificates for Convex Minimization with Inexact Oracle}


\author{Egor Gladin \and Alexander Gasnikov \and Pavel Dvurechensky}

\institute{E. Gladin, Corresponding author \at
              Humboldt-Universit\"at zu Berlin \\
              Unter den Linden 6, 10117 Berlin, Germany \\
              egor.gladin@student.hu-berlin.de
           \and
           A. Gasnikov \at 
Moscow Institute of Physics and Technology \\
9 Institutskiy per., Dolgoprudny, Moscow Region, 141701, Russian Federation \\
Skolkovo Institute of Science and Technology \\
Bolshoy Boulevard 30, bld. 1, 121205
Moscow, Russia \\
Institute for Information Transmission Problems RAS \\
Bolshoy Karetny per. 19, build.1, 127051 Moscow, Russian Federation \\
gasnikov@yandex.ru
\and
             P. Dvurechensky  \at
             Weierstrass Institute for Applied Analysis and Stochastics\\
             Mohrenstr. 39, 10117 Berlin, Germany \\            pavel.dvurechensky@wias-berlin.de
}

\date{}

\maketitle

\begin{abstract}

    Accuracy certificates for convex minimization problems allow for online verification of the accuracy of approximate solutions and provide a theoretically valid online stopping criterion. When solving the Lagrange dual problem, accuracy certificates produce a simple way to recover an approximate primal solution and estimate its accuracy. In this paper, we generalize accuracy certificates for the setting of inexact first-order oracle, including the setting of primal and Lagrange dual pair of problems. We further propose an explicit way to construct accuracy certificates for a large class of cutting plane methods based on polytopes. As a by-product, we show that the considered cutting plane methods can be efficiently used with a noisy oracle even thought they were originally designed to be equipped with an exact oracle. Finally, we illustrate the work of the proposed certificates in the numerical experiments highlighting that our certificates provide a tight upper bound on the objective residual. 

\end{abstract}


\section{Introduction}

The authors of \cite{nemirovski2010accuracy} introduced the notion of accuracy certificates for convex minimization and other problems with a convex structure. These certificates verify the accuracy of an approximate solution at any stage of an optimization algorithm execution. 
Although many algorithms have convergence rate estimates, those often involve parameters unknown in practice, e.g., a constant of Lipschitz continuity of the objective, distance from the starting point to the closest solution, and so on. Accuracy certificates, on the contrary, verify the accuracy of an approximate solution without additional a priori information about the particular problem. Moreover, the accuracy can be verified online and on the fly using the already available information generated by the algorithm.
Thus, the accuracy certificates provide a theoretically valid and practical stopping criterion.
Furthermore, certificates allow an external recipient to verify the accuracy guarantees, without knowing how the algorithm works. This can be useful in some cases where privacy is a priority.

Certificates are extremely useful when an algorithm is applied to the Lagrange dual optimization problem. In this case, they can be used to convert an $\epsilon$-optimal dual solution into an $\epsilon$-optimal solution to the primal problem. Moreover, this approach allows to reuse the information already generated by the algorithm and the approximate primal solution is reconstructed in a direct way, without the knowledge of additional problem parameters.
Remarkably, certificate-based approach allows one to circumvent the following disadvantages of the approach based on the regularization of the primal problem \cite{devolder2012double,gasnikov2016efficient}.
\begin{itemize}
    \item The regularization approach uses an upper bound on the norm of a primal solution. In many cases, such a bound is not available or overestimated which leads to slower convergence.
    \item It requires the target accuracy to be fixed in advance. This raises difficulties when the time limit is exceeded before this accuracy is reached, or when the user decides in the middle of the process that a higher target accuracy is needed.
    \item In some cases, regularization is not enough to reconstruct the primal solution. For example, one needs to impose a $\beta$-ergodicity assumption when applying a dual approach to optimizing a constrained Markov decision process \cite{gladin2023algorithm}.
    \item In some cases, the accuracy deteriorates when one converts a dual point to a primal solution in this way. For example, accuracy $\epsilon$ of a dual solution might result in accuracy $\sqrt{\epsilon}$ of the respective primal solution \cite{gladin2023algorithm}.
\end{itemize}

\textbf{Related literature.}
The paper \cite{nemirovski2010accuracy}
provides a way to compute accuracy certificates for the ellipsoid method. After $\tau$ iterations of the method applied to an $n$-dimensional problem, computation of certificates requires $O(n^3)+O(\tau n^2)$ arithmetic operations (a.o.).
Moreover, the authors mention that one can compute certificates for other cutting-plane methods in a similar fashion by approximating localizer sets with John ellipsoids (for a notion of John ellipsoid see, e.g., \cite{boyd2004convex}, Chapter 8.4). However, we didn't find any uses of this procedure with algorithms other than the ellipsoid method. A possible reason for this is the high computational cost of approximating John ellipsoids \cite{khachiyan1990complexity,nemirovski1999self,anstreicher2002improved,kumar2005minimum,todd2007khachiyan,cohen2019near}.
Fortunately, we show that there is a way to build certificates for polytope-based cutting plane algorithms in a straightforward way by approximately solving a single linear problem, which according to \cite{van2020deterministic} takes only $\widetilde{O}\left(n^{\omega} \right)$ a.o., where  $\widetilde{O}$ hides polylog$(n)$ factors.
We also note that the certificates proposed in \cite{nemirovski2010accuracy} are constructed under the assumption that the first-order information, i.e., subgradients, in the problem are available exactly, which may not always happed in practice.

\textbf{Contributions} of this paper are as follows:
\begin{itemize}
    \item Generalizing the work \cite{nemirovski2010accuracy}, we investigate the properties of accuracy certificates in the setting of minimization problems with \textit{inexact first-order oracle};
    \item We develop a simple and efficient way to obtain accuracy certificates for a large class of cutting plane methods, including Vaidya's method \cite{vaidya1989new,vaidya1996new}, Atkinson-Vaidya algorithm \cite{atkinson1995cutting} and many others;
    \item We show that the considered methods can be efficiently used with a noisy oracle even if they were originally designed to be used with an exact oracle;
    \item  We consider convex problems with (possibly nonlinear) convex inequality constraints and establish a straightforward way to obtain an approximate primal solution based on the information obtained by a method with certificates applied to the dual problem. Generalizing \cite{nemirovski2010accuracy}, we consider nonlinear constraints and allow for inexact solutions of auxiliary problems in each iteration.
\end{itemize}

The rest of the paper is organized as follows. In Section \ref{S:certificates}, we state the minimization problem, and define the separation oracle and the inexact first-order oracle that are used in the algorithms. We also formally define certificates and prove their main property, namely, an upper bound for the objective residual based on a certificate. Section \ref{S:primal-dual} is devoted to the primal-dual setting where we consider convex optimization problems with convex inequality constraints. In particular, we construct separation and inexact oracles for this setting and propose a way to reconstruct an approximate primal solution based on a certificate for the dual problem, with the same accuracy both in terms of the primal objective and constraint violation. In Section \ref{S:certificates_construction}, we describe a wide class of cutting plane methods and propose a way to construct accuracy certificates for these methods. Finally, in Section \ref{S:experiments}, we illustrate the practical efficiency of the proposed certificates.  

\section{Certificates and Their Properties}
\label{S:certificates}

\subsection{Problem Formulation}

Consider a convex minimization problem (CMP)
\begin{equation}\label{eq:cmp}
    \Opt = \min_{x\in X} F(x),
\end{equation}
where
\begin{itemize}
    \item $X \subset \R^n$ is a solid (convex compact set with a nonempty interior) represented by a \textit{Separation oracle} -- a black box which, given on input a point $x \in \R^n$, reports whether or not $x \in \int X$, and in the case of $x \notin \int X$, returns a \textit{separator} -- a vector $e \neq 0$ such that $\langle e, y-x\rangle \leq 0$ for all $y \in X$.
    \item $F: X \rightarrow \R \cup\{+\infty\}$ is a convex function with $\operatorname{Dom}(F)=\{x: F(x)<\infty\} \supseteq \int X$; this function is represented by $\delta$-\textit{oracle} -- a black box which, given on input a point $x \in \operatorname{int} X$, returns a value $\tilde{F}(x)$ such that $|\tilde{F}(x)-F(x)|<\delta$, and a $\delta$-subgradient $\tilde{F}^{\prime}(x)\in \partial_\delta F(x)$ of $F$ at $x$, i.e., a vector $\tilde{F}^{\prime}(x)$ satisfying
    \begin{equation}\label{eq:delta_subgrad}
        F(y) \geq F(x) + \la \tilde{F}^{\prime}(x), y-x\ra -\delta\quad \forall y \in X.
    \end{equation}
\end{itemize}
A point $x \in \operatorname{int} X$ is called a \textit{strictly feasible} solution to \eqref{eq:cmp}. A proximity measure for such a point $x$ to optimality is defined by
$$
\epsilon_{\mathrm{opt}}(x)=F(x)-\inf _{y \in X} F(y)=F(x)-\text { Opt. }
$$
A strictly feasible point $x$ is called $\epsilon$-optimal for \eqref{eq:cmp}, if $\epsilon_{\mathrm{opt}}(x) \leq \epsilon$, that is, if $F(x) \leq \mathrm{Opt}+\epsilon$.

\subsection{Certificates for Convex Minimization Problems}

A computational method for solving the problem \eqref{eq:cmp} within a prescribed accuracy $\epsilon>0$ produces execution protocols $P_\tau=\left\{\left(x_t, e_t\right)\right\}_{t=1}^\tau$, where
\begin{itemize}
    \item $\tau \in \mathbb{N}$ is the current number of steps, 
    \item $x_t \in \R^n$ are the search points generated so far, 
    \item $e_t$ is either a nonzero vector, reported by the Separation oracle and separating $x_t$ and $X$ (this is the case at a \textit{nonproductive} steps $t$ -- those with $\left.x_t \notin \operatorname{int} X\right)$, or is a $\delta$-subgradient $\tilde{F}^{\prime}\left(x_t\right)$ of $F$ at $x_t$ reported by the $\delta$-oracle (this is the case at \textit{productive} steps $t$ -- those with $x_t \in \operatorname{int} X$).
\end{itemize}
The range $1 \leq t \leq \tau$ of the values of $t$ associated with an execution protocol $P_\tau$ is split into the sets $I_\tau, J_\tau$ of indices of productive, resp., nonproductive steps, and the protocol is augmented by the approximate values $\tilde{F}\left(x_t\right)$ of the objective at productive search points $x_t$ -- those with $t \in I_\tau$.
We are about to demonstrate a natural way to certify $\epsilon$-optimality of a strictly feasible solution offered by certificates which are defined as follows.
\begin{definition}
    Let $P_\tau$ be an execution protocol. A \textit{certificate} for this protocol is a collection $\xi=\left\{\xi_t\right\}_{t=1}^\tau$ of weights such that
    \begin{itemize}
        \item $\xi_t \geq 0$ for each $t=1, \ldots, \tau$,
        \item $\sum_{t \in I_\tau} \xi_t=1$.
    \end{itemize}
\end{definition}
Note that certificates exist only for protocols with nonempty sets $I_\tau$.
\begin{definition}
    Given a solid $\mathbf{B}$ known to contain $X$, an execution protocol $P_\tau$ and a certificate $\xi$ for this protocol, we define the quantity
    $$
    \epsilon_{\mathrm{cert}}\left(\xi \mid P_\tau, \mathbf{B}\right) \equiv \max _{x \in \mathbf{B}} \sum_{t=1}^\tau \xi_t\left\langle e_t, x_t-x\right\rangle
    $$
    which we call the \textit{residual of} the certificate $\xi$ on $\mathbf{B}$. Moreover, we define \textit{the approximate solution induced by} $\xi$
    \begin{equation*}
        x^\tau[\xi] := \sum_{t \in I_\tau} \xi_t x_t
    \end{equation*}
    which clearly is a strictly feasible solution to \eqref{eq:cmp}.
\end{definition}
The role of the just defined quantities in certifying accuracy of approximate solutions to \eqref{eq:cmp}
stems from the following
\begin{proposition}\label{prop:resid_def}
    Let $P_\tau$ be a $\tau$-point execution protocol associated with the CMP \eqref{eq:cmp}, $\xi$ be a certificate for $P_\tau$ and $\mathbf{B} \supset X$ be a solid.
    Then
    $x^\tau=x^\tau[\xi]$ is a strictly feasible solution of the
    given CMP, with
    $$		 \epsilon_{\mathrm{opt}}\left(x^\tau\right) \leq \epsilon_{\mathrm{cert}}\left(\xi \mid P_\tau, \mathbf{B}\right) + \delta.
    $$
\end{proposition}
Proof can be found in Appendix \ref{proof:resid_def}.

\section{Recovering Approximate Primal Solution from Dual}
\label{S:primal-dual}

Consider a convex optimization problem with constraints
\begin{equation}\label{eq:primal}
    \Opt = \min_{u \in U} \{f(u) : g(u) \leq 0\},
\end{equation}
where $g(u)=[g_1(u),\ldots, g_n(u)]^\top,\, g_j(u)$ are convex functions, $U$ is a closed convex set.
We assume the problem to be bounded below. A natural way to solve it is to consider its Lagrange dual problem:
\begin{equation}\label{eq:dual}
    \min_{x \geq 0} F(x), \quad F(x)=-\min_{u \in U}\{\underbrace{f(u)+\la x, g(u)\ra}_{\phi(u,x)} \}.
\end{equation}
Assuming that \eqref{eq:primal} satisfies the Slater condition (so that \eqref{eq:dual} is solvable) and that we have at our disposal an upper bound $L$ on the norm $\left\|x_*\right\|_p$ of an optimal solution $x_*$ to \eqref{eq:dual}, we can reduce the problem to solving the following CMP:
\begin{equation}\label{eq:dual2}
    \min_{x \in X} F(x), \quad X = \left\{x \geq 0:\|x\|_p \leq L+1\right\}. 
\end{equation}
We further assume that $\phi(\cdot,x)$ is bounded from below for every $x\in X$, i.e., $X \subseteq \dom F \equiv \{x\in \R^n: \min_{u \in U}\phi(u,x) > -\infty \}$.
This is the case, for example,
if $U$ is compact or if $f(u)$ is strongly convex.

\subsection{Separation and $\delta$-oracles}

Separation oracle for $X$ is easily constructed in the following way: let $x'\notin \int X$. If $x_i' \leq 0$, then the vector $-e_i$ (having $-1$ is position $i$ and 0 in others) is a separator since for any $x\in X$ it holds $-e_i^\top x \leq 0 \leq -x_i' = -e_i^\top x'$. If $\|x'\|_p \geq L+1$, then let $a\in \R^n$ be the vector satisfying
$$
\operatorname{sign} a_i = \operatorname{sign} x_i',\; |a_i|^q=\textstyle \frac{|x_i'|^p}{\|x'\|_p^p},\; i=1,\ldots, n,
$$
so that Hölder's inequality for $a$ and $x$ becomes an equality: $a^\top x' = \|a\|_q \|x'\|_p \geq L+1$ since $\|a\|_q=1$. Thus, $a$ is a separator since for any $x\in X$ it holds $a^\top x \leq \|a\|_q \|x'\|_p \leq L+1 \leq a^\top x'$.

Let us now show that it is easy to equip $F$ with a $\delta$-oracle provided that the aforementioned assumptions hold and that an efficient first-order method for solving the convex problem $\min_{u \in U}\phi(u,x)$ up to a prescribed accuracy $\delta$ is available.
Let $u_x$ be the point returned by such method,
i.e., $\phi(u_x,x) - F(x) \leq \delta$ (we also write: $u_x \in \arg\underset{u\in U}{\min}^\delta\phi(u,x)$). It follows from the argument on page 132 of \cite{polyak1983intro} that $-g(u_x)\in \partial_\delta F(x)$. We provide this argument below:
\begin{align*}
    \forall x'\in X,\quad F(x')&=-\min_{u \in U}\phi(u,x')\geq -\phi(u_x,x') = -\phi(u_x,x)-\la g(u_x), x'-x\ra \\
    &\geq F(x)-\la g(u_x), x'-x\ra -\delta.
\end{align*}

\subsection{Reconstructing Primal Solution}

With separation and $\delta$-oracles at hand, we can solve the dual problem \eqref{eq:dual2}. It turns out that accuracy certificates allow us to recover nearly feasible and nearly optimal solution for \eqref{eq:primal}. The following statement generalizes Proposition 5.1 from \cite{nemirovski2010accuracy}.
\begin{proposition}\label{prop:primal_dual}
    Let \eqref{eq:dual2} be solved by a black-box-oriented method, $P_\tau=\left\{I_\tau, J_\tau,\left\{x_t, e_t\right\}_{t=1}^\tau\right\}$ be the execution protocol upon termination, with
    \begin{equation*}
        e_t = -g(u_t),\quad u_t \in \arg\underset{u\in U}{\min}^\delta\phi(u,x^t),\quad t \in I_\tau.
    \end{equation*}
    Let also $\xi$ be an accuracy certificate for this protocol. Set
    $\hat{u}=\sum_{t \in I_\tau} \xi_t u_{t}$,
    then $\hat{u} \in U$ and
    \begin{align}
        \left\|[g(\hat{u})]_{+}\right\|_q &\leq \epsilon_{\mathrm{cert}}\left(\xi \mid P_\tau, X\right) + \delta, \label{eq:constr_viol} \\
        f(\hat{u})-\Opt &\leq \epsilon_{\mathrm{cert}}\left(\xi \mid P_\tau, X\right) + \delta, \label{eq:opt_gap}
    \end{align}
    where $[g(\hat{u})]_{+}$ is the ``vector of constraint violations'' obtained from $g(\hat{u})$ by replacing the negative components with 0, and $q=p /(p+1)$.
\end{proposition}
Proof can be found in Appendix \ref{proof:primal_dual}.

Proposition \ref{prop:primal_dual} shows that the vector $\hat{u}=\sum_{t \in I_\tau} \xi_t u_{t}$
is nearly feasible and nearly optimal for \eqref{eq:primal}, provided that $\epsilon_{\mathrm{cert}}\left(\xi \mid P_\tau, \mathbf{B}\right)$ is small.

\section{Accuracy Certificates for Cutting Plane Methods}
\label{S:certificates_construction}

\subsection{Generic Cutting Plane Algorithm with $\delta$-Oracle}

A generic cutting plane algorithm with $\delta$-oracle, as applied to a CMP \eqref{eq:cmp}, 	
builds a sequence of search points $x_t \in \R^n$ along with a sequence of localizers $Q_t$ -- solids such that $x_t \in \operatorname{int} Q_t, t=1,2, \ldots$. The algorithm is as follows:

\textbf{Initialization:} Choose a solid $Q_1 \supset X$ and a point $x_1 \in \operatorname{int} Q_1$.

\textbf{Step} $t=1,2, \ldots$: given $x_t, Q_t$,
\begin{enumerate}
    \item Call Separation oracle, $x_t$ being the input. If the oracle reports that $x_t \in \operatorname{int} X$ (productive step), go to 2. Otherwise (nonproductive step) the oracle reports a separator $e_t \neq 0$ such that $\left\langle e_t, x-x_t\right\rangle \leq 0$ for all $x \in X$. Go to 3.
    \item Call $\delta$-oracle to compute $e_t=\tilde{F}^{\prime}(x_t)\in \partial_\delta F(x_t)$. If $e_t=0$, terminate, otherwise go to 3.
    \item \label{step_3} Set
    
    $$\widehat{Q}_{t+1}=\left\{x \in Q_t:\left\langle e_t, x-x_t\right\rangle \leq 0\right\} .$$
    
    Choose as $Q_{t+1}$, a solid which contains the solid $\widehat{Q}_{t+1}$. Choose $x_{t+1} \in \operatorname{int} Q_{t+1}$ and loop to step $t+1$.
\end{enumerate}
For a solid $B \subset \R^n$, let $\rho(B)$ be the radius of Euclidean ball in $\R^n$ with the same $n$-dimensional volume as the one of $B$. A cutting plane algorithm with $\delta$-oracle applied to the problem \eqref{eq:cmp} is called converging if for the associated localizers $Q_t$ one has $\rho\left(Q_t\right) \rightarrow 0,\, t \rightarrow \infty$. Some examples of converging cutting plane algorithms are the center of gravity method \cite{levin1965minimization,newman1965location}, the ellipsoid method \cite{yudin1976informational,shor1977cutting}, the inscribed ellipsoid algorithm \cite{khachiyan1988method}, the circumscribed simplex algorithm \cite{bulatov1982method,yamnitsky1982old}, Vaidya’s algorithm  \cite{vaidya1989new,vaidya1996new}.

\subsection{Polytope-Based Cutting Plane Algorithms}\label{subsec:polytope_algo}
Recall that a full-dimensional polytope is a bounded set with nonempty interior of the form
$$
    Q(A, b) = \{x\in\R^n: a_i^\top x \leq b_i, i=1,\ldots,m\} = \{x\in\R^n: A x\leq b\}
$$
for given
$$
    A=\left[\begin{array}{c} a_1^\top \\ \vdots \\ a_{m}^\top \end{array}\right],\; a_1, \ldots, a_{m}\in\R^n,\; b\in\R^{m}.
$$
In what follows, we consider implementations of a generic cutting plane algorithm with $\delta$-oracle where localizers are full-dimensional polytopes, i.e., $Q_t = Q(A_t, b_t)$. In what follows, we omit the subscript $t$ for brevity when it doesn't cause ambiguity, i.e., we write $Q_t=Q(A,b)$ and implicitly assume that $m, A$ and $b$ depend on $t$.

Moreover, we assume that if the constraint $a_i^\top x \leq b_i$ was added at the step $t(i)$, then $a_i = e_{t(i)}$. Note that item \ref{step_3} in the description of the generic cutting plane algorithm with $\delta$-oracle implies that $a_i^\top x_{t(i)} \leq b_i$, i.e.,
when a new constraint is added, the current iterate satisfies it. We will refer to the group of methods described above as \textit{polytope-based cutting plane algorithms with $\delta$-oracle}.

\subsection{Building Accuracy Certificates}

Consider a \textit{nonterminal} step $\tau$ (i.e., the one with $e_\tau \neq 0$) of a polytope-based cutting plane algorithms with $\delta$-oracle. The respective localizer $Q_{\tau+1}$ is formed by the set of constraints $a_i^\top x \leq b_i, i=1,\ldots,m$ which can be divided into three disjoint sets:
$\{1, \ldots, m\} = \It \cup \Pt \cup \Nt$, where
\begin{itemize}
    \item $\It$ (not to be confused with $I_\tau$) corresponds to \textbf{I}nitial constraints that were present in $Q_1$,
    \item $\Pt$ (not to be confused with $P_\tau$) corresponds to constraints added during \textbf{P}roductive steps of the algorithm,
    \item $\Nt$ corresponds to constraints added during \textbf{N}onproductive steps.
\end{itemize}
Note that if a constraint was removed during execution of the algorithm, it does not appear in any of the sets $\It, \Pt, \Nt$.

The following LP problem will play a crucial role in building certificates:
\begin{align}\label{eq:lp}
    \max_{\lambda \in \R^{m}}\, D_{\tau}(\lambda) &:= \sum_{i\in \Pt} \lambda_i \| a_i \|_2, \\
    \text{s.t. } \lambda &\geq 0,\nonumber \\
    A^\top \lambda &= 0,\nonumber \\
    b^\top \lambda &\in [0,2].\nonumber
\end{align}
\begin{lemma}\label{lem:boundedness}
    The LP problem \eqref{eq:lp} is feasible and bounded.
\end{lemma}
Proof can be found in Appendix \ref{proof:boundedness_lemma}.
\begin{definition}\label{def:cert}
    If $\lambda$ is a feasible point in the LP problem \eqref{eq:lp} and $d_\tau:= \sum_{i\in \Pt} \lambda_i >0$, define $\xi=\{\xi_t\}_1^\tau$ as follows:
    \begin{enumerate}
        \item For every $i\in \Pt\cup \Nt$, set $\xi_{t(i)}:=\frac{\lambda_i}{d_\tau}$, where $t(i)$ is the step when the constraint $a_i^\top x\leq b_i$ was added,
        \item For all other steps $t$, set $\xi_{t}:=0$
    \end{enumerate}
\end{definition}
Observe that this definition implies
\begin{equation}\label{eq:D_tau}
    D_\tau(\lambda) = \sum_{i\in \Pt} \lambda_i \| a_i \|_2 = d_\tau \sum_{i\in \Pt} \xi_{t(i)} \| e_{t(i)} \|_2 = d_\tau \sum_{t\in I_\tau} \xi_{t} \| e_{t} \|_2.
\end{equation}
In what follows, we sometimes write $D_\tau$ in place of $D_\tau(\lambda)$ for brevity.
\begin{lemma}\label{lem:resid}
    If $\lambda$ is a feasible point for the LP problem \eqref{eq:lp} with $d_\tau>0$, then $\xi$ from Definition \ref{def:cert} is a certificate.
    If $\epsilon_\tau := \frac{2}{D_\tau}<r$, then
    \begin{equation}\label{eq:resid2}
        \resid\resar \leq \frac{\epsilon_\tau}{r-\epsilon_\tau} W_\tau,
    \end{equation}
    where
    \begin{equation}\label{eq:variation}
        W_\tau :=\max _{t \in I_\tau} \max _{x \in X}\left\langle e_t, x-x_t\right\rangle,
    \end{equation}
    and $r=r(X)$ is the largest of the radii of Euclidean balls contained in $X$.
\end{lemma}
Proof can be found in Appendix \ref{proof:resid_lemma}.
\begin{remark}
    Informally speaking, inequality \eqref{eq:resid2} shows that the larger $D_\tau$ is, the more accurate the estimate $\resid\resar$ is, provided that $W_\tau$ is bounded.
\end{remark}

\begin{theorem}\label{thm:D_lower}
    An optimal solution $\lambda^*$ for the LP problem \eqref{eq:lp} satisfies
    \begin{equation}\label{eq:D_lower}
        D_\tau(\lambda^*) \geq D^{-1}\left(Q_1\right)\left(\frac{r}{2 n \rho\left(Q_{\tau+1}\right)}-1\right),
    \end{equation}
    where $D\left(Q_1\right)$ is the Euclidean diameter of $Q_1$.
\end{theorem}
Proof can be found in Appendix \ref{proof:D_lower}.
\begin{remark}
    Since the quantity $D_\tau$ is always nonnegative, the inequality \eqref{eq:D_lower} can only be useful when $\rho\left(Q_{\tau+1}\right) < \frac{r}{2n}$.
\end{remark}
Before we move on to the most important corollary, let us mention that the convergence rate of a cutting plane method is basically described by how fast $\rho\left(Q_{\tau+1}\right)$ is decreasing as $\tau$ grows. For example, for Vaidya's method 
$$
    \rho\left(Q_{\tau+1}\right) \leq C_1\cdot \rho\left(Q_{1}\right) e^{-C_2\tau/n}
$$
for some $C_1, C_2>0$. It can be shown that this estimate implies $\epsilon_{\mathrm{opt}}\left(x^\tau\right) = O\left(e^{-C_2\tau/n}\right)$, where $x^\tau$ is a point returned by Vaidya's method after $\tau$ iterations.
\begin{corollary}\label{cor:convergence}
    Let $\alpha\in [0, 1)$ be the relative accuracy in the LP problem \eqref{eq:lp}.
    If $\tau$ is a nonterminal iteration number of a polytope-based cutting plane algorithms with $\delta$-oracle such that
    $$
    \rho\left(Q_{\tau+1}\right) \leq \frac{(1-\alpha) r^2}{16 n D\left(Q_1\right)},
    $$
    and $\lambda$ is a feasible point for the LP problem \eqref{eq:lp} with $D_\tau(\lambda) \geq (1-\alpha) D_\tau(\lambda^*)$, then the respective certificate $\xi$ is well defined, and
    $$
    \resid\resar \leq \frac{16 n D\left(Q_1\right) W_\tau}{(1-\alpha) r^2} \rho\left(Q_{\tau+1}\right).
    $$
    In particular, if $\sup_{x,y\in X}\left(F(x) - F(y) \right) \leq C < \infty$, then $W_\tau \leq C+\delta$ and
    $$
    \resid\resar \leq \frac{16 n D\left(Q_1\right) (C+\delta)}{(1-\alpha) r^2} \rho\left(Q_{\tau+1}\right).
    $$
\end{corollary}
Proof can be found in Appendix \ref{proof:convergence}.
\begin{remark}
    Parameter $\alpha\in (0, 1]$ provides a trade-off between the number of iterations performed by a cutting plane method and the accuracy of solving the LP problem \eqref{eq:lp}.
\end{remark}
\begin{remark}
    Complexity of constructing certificates is, in essence, the complexity of solving LP \eqref{eq:lp} up to a chosen relative accuracy $\alpha$ (say, $\alpha=1/2$).
    When a method uses polytopes formed by $\widetilde{O}(n)$ of constraints (which is the case, for example, for Vaidya's method), the LP \eqref{eq:lp} can be solved in $\widetilde{O}\left(n^\omega \log (n / \alpha)\right)$ arithmetic operations \cite{van2020deterministic}.
    Here $\widetilde{O}$ hides polylog$(n)$ factors, $O\left(n^\omega\right)$ is the time required to multiply two $n \times n$ matrices.
\end{remark}
\begin{remark}
    Corollary \ref{cor:convergence} implies that all polytope-based cutting plane algorithms
    can be used with a $\delta$-oracle, $\delta \in [0, \varepsilon)$, to achieve an $\varepsilon$-optimal solution provided that their localizers' volumes converge to zero.
\end{remark}

\section{Numerical Experiments}
\label{S:experiments}

We present the results of numerical experiments which aim to showcase the performance of certificates described in the previous section and compare it to that of the certificates from the paper \cite{nemirovski2010accuracy}.
Vaidya's cutting plane method \cite{vaidya1989new,vaidya1996new} is chosen to demonstrate the certificates from Definition \ref{def:cert} in action. Such a choice is made because it fulfills requirements of Subsection \ref{subsec:polytope_algo}, in particular, its localizers are polytopes. Moreover, it is the first optimal cutting plane method in terms of the number of oracle calls.
The ellipsoid method is used to show the performance of the certificates based on Algorithm 4.2 from \cite{nemirovski2010accuracy} which was designed for this method. Although authors mention that it is possible to adapt the certificate computation procedure to other methods, details are omitted. Furthermore, we didn't find any uses of this procedure with algorithms other than the ellipsoid method.

Consider the following nonsmooth convex optimization problem taken from the book \cite{nesterov2018lectures} (subsection 3.2.1):
\begin{equation}\label{eq:experim}
    \min_{x\in \mathbb{R}^n} \Bigl\{ F(x):= \max_{i=1,\ldots,n}x_i + \frac{\mu}{2}\|x\|_2^2 \Bigr\}.
\end{equation}
As proposed in the book, we take the initial point to be $x^0=0$ and let the first-order oracle called at a point $x$ return (apart from the function value) a subgradient $f'(x) = e_{i_*} + \mu x$, where $i_*:=\min\{j\, |\, x_j= \displaystyle\max_{i=1,\ldots,n}x_i \}$. Note that the problem has a closed-form solution $x_* = -\frac{1}{\mu n}\mathbf{1}$ (see \cite{nesterov2018lectures}), which makes it possible to compute the quantities $\epsilon_{\mathrm{opt}}(x^\tau)$ in the experiment. We turn \eqref{eq:experim} into a problem on a solid by setting $X$ to be a Euclidean ball of radius $10\cdot \|x_*\|_2$ centered at the origin.

\begin{figure}[h]
\label{fig:exp}
\centering
\includegraphics[width=\textwidth]{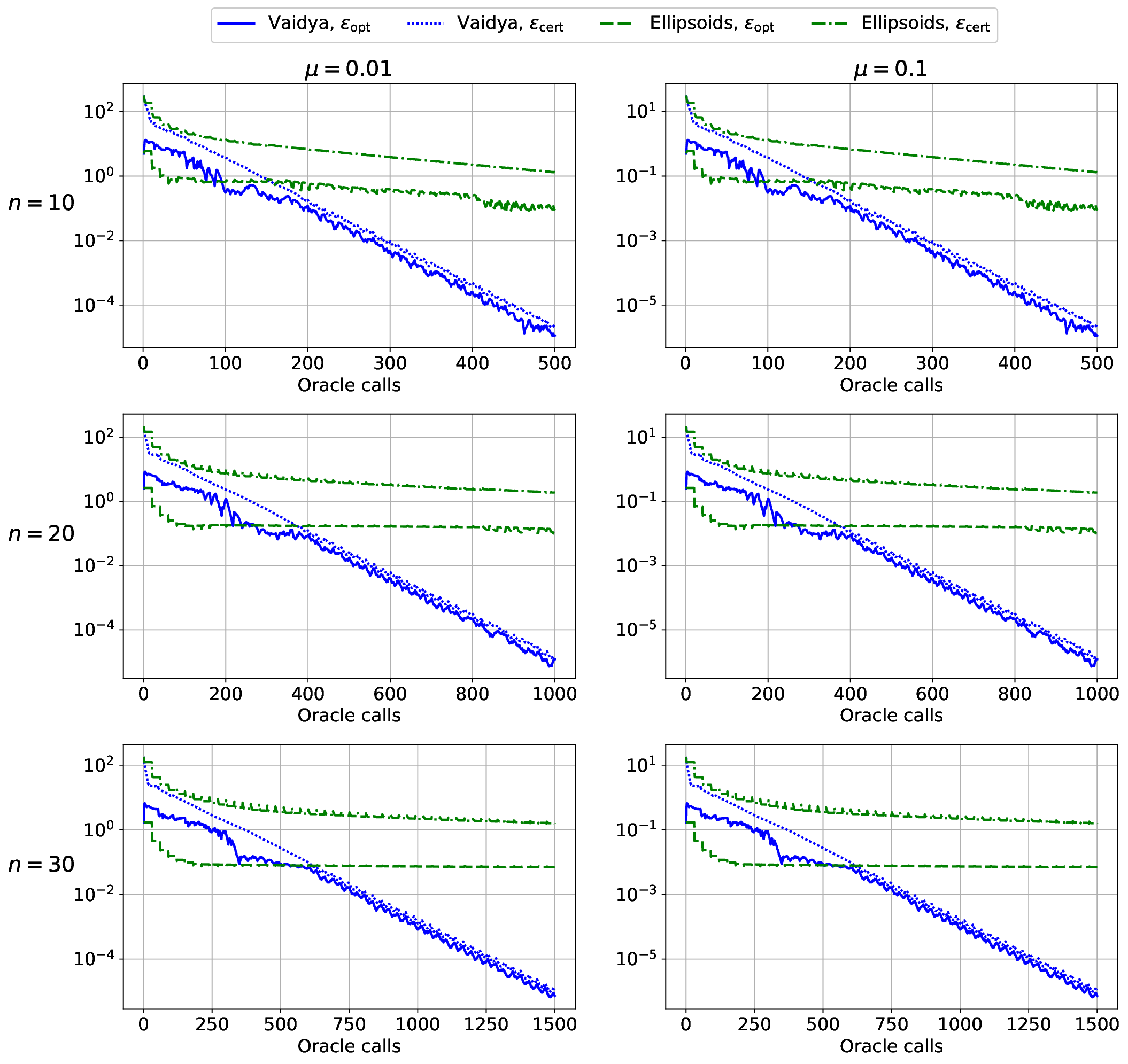}
\caption{Performance of Vaidya's and the ellipsoid methods for problem \eqref{eq:experim}. The rows correspond to different dimensions of the problem ($n=10,20,30$). The left and right columns present the results with small ($\mu=0.01$) and medium ($\mu =0.1$) regularization, respectively. X-axis represents number of oracle calls. Solid and dashed lines depict $\epsilon_{\mathrm{opt}}$ for Vaidya's and the ellipsoid methods, respectively. Dotted and dash-dotted lines depict $\resid$ for Vaidya's and the ellipsoid methods, respectively.}
\end{figure}

Figure \ref{fig:exp} presents the results of the experiments.
The rows represent different dimensions of the problem \eqref{eq:experim} ($n=10,20,30$). The left and right columns correspond to small ($\mu=0.01$) and medium ($\mu =0.1$) regularization, respectively. X-axis depicts the number of oracle calls. Solid and dashed lines represent $\epsilon_{\mathrm{opt}}$ for Vaidya's and the ellipsoid methods, respectively. Dotted and dash-dotted lines depict $\resid$ for Vaidya's and the ellipsoid methods, respectively.

As we see from Figure \ref{fig:exp}, the certificates from Definition \ref{def:cert} for Vaidya's method provide an upper bound $\resid$ on the optimality gap $\epsilon_{\mathrm{opt}}$ which becomes tight after a few hundred iterations. The certificates based on Algorithm 4.2 from \cite{nemirovski2010accuracy} for the ellipsoid method yield a less tight bound. We refer the reader to Section \ref{sec:concl} for discussion of this phenomenon. Moreover, the figure illustrate the fact that Vaidya's method scales much better with the dimension $n$.

\subsection{Implementation Details}

We use the version of Vaidya's cutting plane method from the paper \cite{anstreicher1997vaidya} since it is more practical than the original version. The parameters used are $\varepsilon=5\cdot 10^{-3}$, $\tau=1$, see the aforementioned paper for details.
Certificates and their residuals were computed after each iteration for illustration purposes.
The experiments were conducted using programming language Python 3.11.5 with packages numpy v1.26.0 and scipy v1.11.3.
The source code is available at \url{https://github.com/egorgladin/vaidya-with-certificates}.

\section{Conclusions}\label{sec:concl}

The present paper generalizes the notion of accuracy certificates to the case of convex optimization problems with inexact oracle and establishes properties of such certificates. In particular, we show how they provide a simple way to recover primal solutions when solving a wide class of Lagrange dual problems. Additionally, we develop a new recipe to construct certificates suitable for cutting plane methods which use polytopes as localizers. A prominent example is Vaidya's method which is asymptotically optimal in terms of the number of oracle calls. Arithmetic complexity of our recipe is equivalent to that of approximately solving an LP problem. Notably, the requirements for the accuracy of such approximate solutions are very mild. As of this writing, this can be done in current matrix multiplication time.
As an important by-product of our analysis, we conclude that all polytope-based cutting plane algorithms can be used with an inexact oracle to achieve a near-optimal solution provided that their localizers' volumes converge to zero.

Numerical experiments show that the proposed procedure for computing certificates may be superior to the existing approach which we build on.
A possible reason for this phenomenon is that we look for certificates that directly maximize a function used to bound the residual. The previous approach, in contrast, simply produces a point with a sufficiently large value of that function.
Although we use such a point in the analysis to bound the optimal value from below, the maximum may turn out to be considerably larger, which leads to better residuals.

A possible way to improve the presented approach for computing certificates is to introduce warm starts. Namely, one could solve the LP problem once, and then update the solution after each iteration as the problem is slightly modified between two consecutive iterations. This may further reduce the complexity of the approach.
Another important direction of future research is the exploration of accuracy certificates for variational inequalities with inexact oracle.

\begin{acknowledgements}
The work of E. Gladin and P. Dvurechensky is funded by the Deutsche Forschungsgemeinschaft (DFG, German Research Foundation) under Germany's Excellence Strategy – The Berlin Mathematics Research Center MATH+ (EXC-2046/1, project ID: 390685689).

The work of A. Gasnikov was supported by Ministry of Science and Higher Education grant No. 075-10-2021-068.
\end{acknowledgements}

\appendix  
\section{Proofs of Propositions and Lemmas}


\subsection{Proof of Proposition \ref{prop:resid_def}}\label{proof:resid_def}
Since the points $x_t$ with $t \in I_t$ belong to $\operatorname{int} X$ and $X$ is convex, $x^\tau$ (which is a convex combination of these points) belongs to $\operatorname{int} X$ and thus is a strictly feasible solution.

Define 
\begin{align}
    F_*\left(\xi \mid P_\tau, \mathbf{B}\right) & \equiv \min _{x \in \mathbf{B}}\left[\sum_{t \in I_\tau} \xi_t\left[F\left(x_t\right)+\left\langle e_t, x-x_t\right\rangle\right]+\sum_{t \in J_\tau} \xi_t\left\langle e_t, x-x_t\right\rangle\right] \nonumber \\
    & =\sum_{t \in I_\tau} \xi_t F\left(x_t\right)-\epsilon_{\mathrm{cert}}\left(\xi \mid P_\tau, \mathbf{B}\right) . \label{eq:def_Fstar}
\end{align}
We will first show that
\begin{equation}\label{eq:F_star}
    F_*\left(\xi \mid P_\tau, \mathbf{B}\right) \leq \text{Opt} + \delta.
\end{equation}
Let $x \in X$. Then, due to the origin of vectors $e_t$, we have $\left\langle e_t, x-x_t\right\rangle \leq 0$ for $t \in J_\tau$ and $F\left(x_t\right)+\left\langle e_t, x-x_t\right\rangle \leq F(x) + \delta$ for $t \in I_\tau$. Taking weighted sum of these inequalities with the weights determined by a certificate $\xi$, we get
$$
\sum_{t \in I_\tau} \xi_t\left[F\left(x_t\right)+\left\langle e_t, x-x_t\right\rangle\right]+\sum_{t \in J_\tau} \xi_t\left\langle e_t, x-x_t\right\rangle \leq F(x) + \delta.
$$
Hence, taking the infimum of both sides over $x \in X \cap \operatorname{Dom} F$,
$$
\min _{x \in X}\left[\sum_{t \in I_\tau} \xi_t\left[F\left(x_t\right)+\left\langle e_t, x-x_t\right\rangle\right]+\sum_{t \in J_\tau} \xi_t\left\langle e_t, x-x_t\right\rangle\right] \leq \mathrm{Opt} + \delta.
$$
It remains to note that the left hand side in this inequality is $\geq F_*\left(\xi \mid P_\tau, \mathbf{B}\right)$ due to $X \subseteq \mathbf{B}$.

Now, observe that
\begin{align*}
    \epsilon_{\mathrm{opt}}\left(x^\tau\right) &= F(x^\tau)-\mathrm{Opt} \stackrel{\text{convexity}}{\leq} \sum_{t \in I_\tau}\xi_t F\left(x_t\right) -\mathrm{Opt} \\
    &\stackrel{\eqref{eq:F_star}}{\leq} \sum_{t \in I_\tau}\xi_t F\left(x_t\right) - F_*\left(\xi \mid P_\tau, \mathbf{B}\right) + \delta \\
    &\stackrel{\eqref{eq:def_Fstar}}{=} \epsilon_{\mathrm{cert}}\left(\xi \mid P_\tau, \mathbf{B}\right) + \delta.
\end{align*}

\subsection{Proof of Proposition \ref{prop:primal_dual}}\label{proof:primal_dual}

First, let us provide a lower bound on the certificate residual.
\begin{align*}
    \epsilon_{\mathrm{cert}}&\left(\xi \mid P_\tau, X\right)\geqslant && \\
    \geqslant& \max _{x \in X} \sum_{t=1}^\tau \xi_t\left\la e_t, x_t-x\right\ra \\
    \geqslant& \max _{x \in X} \sum_{I_\tau} \xi_t\left\la e_t, x_t-x\right\ra && \Big|\; \left\la e_t, x-x_t\right\ra \leq 0\; \forall t \in J_\tau, x\in X  \\
    =&-\sum_{I_\tau} \xi_t\left\la g\left(u_t\right), x_t\right\ra+\max _{x \in X}\Bigl\la\sum_{I_\tau} \xi_t g\left(u_t\right), x\Bigr\ra && \Big|\; e_t=-g\left(u_t\right)\; \forall t \in I_\tau \\
    \geqslant&-\sum_{I_\tau} \xi_t\left\la g\left(u_t\right), x_t\right\ra+\max _{x \in X}\la g(\hat{u}), x\ra && \Big|\; g_i(\underbrace{\textstyle\sum_{I_\tau} \xi_t u_t}_{\hat{u}}) \stackrel{\text{conv.}}{\leqslant} \sum_{I_\tau} \xi_t g_i\left(u_t\right) \\
    =&-\sum_{I_\tau} \xi_t\left\la g\left(u_t\right), x_t\right\ra + (L+1)\left\|[g(\hat{u})]_{+}\right\|_q. && \Big|\; \max _{x \in X}\la g(\hat{u}), x\ra \stackrel{\text{Hölder}}{=} (L+1)\left\|[g(\hat{u})]_{+}\right\|_q
\end{align*}
Since $u_t \in \arg\underset{u\in U}{\min}^\delta\phi(u,x_t), t\in I_\tau$, it holds
$$
    \phi(u_t,x_t) = f\left(u_t\right)+\left\langle g\left(u_t\right), x_t\right\rangle \leq f(u)+\left\la g(u), x_t\right\ra+\delta,\quad \forall t \in I_\tau,\, u\in U.
$$
Summing over $t\in I_\tau$ and using the inequality $f(\hat{u})\leq \sum_{I_\tau} \xi_t f\left(u_t\right)$, we get
\begin{equation}\label{eq:primal_upper_bound}
    f(\hat{u})-f(u)-\langle g(u), \bar{x}\rangle \leqslant-\sum_{I_\tau} \xi_t\left\langle g\left(u_t\right), x_t\right\rangle+\delta \quad \forall u \in U,
\end{equation}
where $\bar{x}:=\sum_{I_\tau} \xi_t x_t$. Combining the lower bound on the certificate residual and \eqref{eq:primal_upper_bound} where $u=u_*$ (an optimal solution to \eqref{eq:primal}), we arrive at
\begin{align}
    \epsilon_{\mathrm{cert}}\left(\xi \mid P_\tau, X\right) + \delta &\geqslant (L+1)\left\|[g(u)]_{+}\right\|_q+f(\hat{u})-f\left(u_*\right) - \left\langle g\left(u_*\right), \bar{x}\right\rangle \label{eq:primal_dual_bound} \\
    &\geqslant f(\hat{u})-f\left(u_*\right). \nonumber
\end{align}
Thus, \eqref{eq:opt_gap} is established.

Due to Slater's theorem,
\begin{align*}
    f(u_*)=-F(x_*)=\min_{u \in U}\{f(u)+\left\la g(u), x_*\right\ra\} &\leq f(\hat{u})+\left\la g(\hat{u}), x_*\right\ra \\
    &\leq f(\hat{u})+L\left\|[g(\hat{u})]_{+}\right\|_q.
\end{align*}
Combining this inequality with \eqref{eq:primal_dual_bound}, we arrive at \eqref{eq:constr_viol}.

\subsection{Proof of Lemma \ref{lem:boundedness}}\label{proof:boundedness_lemma}

Feasibility is evident since the zero vector satisfies all constraints. Suppose that the feasible set is unbounded, then we will show that there exist a vector $\nu \in \R^m$ such that
\begin{equation}\label{eq:ray}
    \nu\geq 0,\; A^\top\nu=0,\; b^\top\nu = 0,\; \|\nu\|_2=1.
\end{equation}
Indeed, let $\{\lambda_k\}_1^\infty$ be an unbounded sequence of feasible points. Specifically,
$$
    \lambda_k \geq 0,\; A^\top\lambda_k=0,\; b^\top\lambda_k\in [0,2],\; \|\lambda_k\|_2\geq k \quad \forall k \in \mathbb{N}.
$$
Define $\nu_k:=\frac{\lambda_k}{\|\lambda_k\|_2}$, then 
$$
\nu_k \geq 0,\; A^\top\nu_k=0,\; b^\top\nu_k\in \bigl[0,\textstyle\frac{2}{k}\bigr],\; \|\nu_k\|_2=1 \quad \forall k \in \mathbb{N}.
$$
Let $\nu_{k_j}$ be a convergent subsequence, then its limit $\nu$ satisfies \eqref{eq:ray}.

Now, let $Q_{\tau+1}=\{x\in \R^n:\; Ax\leq b\}$ be the current localizer. Since it has nonempty interior, there exist $x_+ \in Q_{\tau+1}$ such that $b-A x_+>0$. At the same time,
\begin{equation}\label{eq:ray_dot}
    \nu^\top(b-Ax_+)=b^\top\nu-x_+^\top A^\top \nu = 0.
\end{equation}
Thus, $\nu=0$ which contradicts $\|\nu\|_2=1$. Therefore, the feasible set is bounded.

\subsection{Proof of Lemma \ref{lem:resid}}\label{proof:resid_lemma}

The fact that $\xi$ is a certificate follows from its construction. Let us first show that $\resid\resar \leq \frac{2}{d_\tau}$. For any $x\in Q_1$, it holds $a_i^\top x \leq b_i,\, \forall i\in \It$. Therefore,
$$
    \lambda^\top(b-Ax) = \sum_{i\in \It\cup \Pt\cup \Nt}\lambda_i (b_i - a_i^\top x) \geq \sum_{i\in \Pt\cup \Nt}\lambda_i (b_i - a_i^\top x).
$$
On the other hand, $\lambda^\top(b-Ax) = \lambda^\top b \leq 2$ since $\lambda$ is a feasible point for LP problem \eqref{eq:lp}. Thus,
\begin{align}
    \sum_{t =1}^\tau \xi_{t} \la e_t, x_t-x\ra &= d_\tau^{-1} \sum_{i\in \Pt\cup \Nt}\lambda_i a_i^\top (x_{t(i)}-x) \nonumber \\
    &\leq d_\tau^{-1} \sum_{i\in \Pt\cup \Nt}\lambda_i (b_i-a_i^\top x) \leq \frac{2}{d_\tau}, \label{eq:almost_resid}
\end{align}
where we used $a_i^\top x_{t(i)} \leq b_i, i\in \Pt\cup \Nt$, see the end of subsection \ref{subsec:polytope_algo}.
We maximize the left-hand side of the last equation w.r.t. $x\in Q_1$ to obtain $\resid\resar \leq \frac{2}{d_\tau}$.

To prove \eqref{eq:resid2}, let $\tau$ be such that $\epsilon := \epsilon_\tau<r$, and let $\bar{x}$ be the center of Euclidean ball $B$ of the radius $r$ which is contained in $X$. Observe that the definition \eqref{eq:variation} of $W_\tau$ implies
$$
t \in I_\tau \Rightarrow\left\langle e_t, x-x_t\right\rangle \leq W_\tau\; \forall x \in B,
$$
hence $\left\langle e_t, \bar{x}-x_t\right\rangle \leq W_\tau-r\left\|e_t\right\|_2$. Recalling what $e_t$ is for $t \in J_\tau \equiv\{1, \ldots, \tau\} \backslash I_\tau$, we get the relations
\begin{equation*}
    \left\langle e_t, \bar{x}-x_t\right\rangle \leq \begin{cases}W_\tau-r\left\|e_t\right\|_2, & t \in I_\tau \\ 0, & t \in J_\tau\end{cases}.
\end{equation*}
In particular,
\begin{equation}\label{eq:ball_center_ineq}
    \la a_i, \bar{x}-x_{t(i)}\ra \leq \begin{cases}W_\tau-r\left\|a_i\right\|_2, & i \in \Pt \\ 0, & t \in \Nt\end{cases}.
\end{equation}
Now let $x \in Q_1$, and let $y=\frac{(r-\epsilon) x+\epsilon \bar{x}}{r}$, so that $y \in Q_1$. By \eqref{eq:almost_resid} we have
$$
\sum_{i\in \Pt\cup \Nt}\lambda_i \la a_i, x_{t(i)}-y\ra \leq 2.
$$
Multiplying this inequality by $r$ and adding weighted sum of inequalities \eqref{eq:ball_center_ineq}, the weights being $\lambda_{i} \epsilon$, we get
$$
\sum_{i\in \Pt\cup \Nt}\lambda_i \langle a_i, \underbrace{r x_{t(i)}-r y+\epsilon \bar{x}-\epsilon x_{t(i)}}_{(r-\epsilon)\left(x_{t(i)}-x\right)}\rangle \leq 2 r+\epsilon W_\tau d_\tau-r \epsilon D_\tau .
$$
The right hand side in this inequality, by the definition of $\epsilon$, is $\epsilon W_\tau d_\tau$, and we arrive at the relation
$$
(r-\epsilon)\cdot \sum_{i\in \Pt\cup \Nt} \lambda_i \left\langle a_i, x_{t(i)}-x\right\rangle \leq \epsilon W_\tau d_\tau \iff \sum_{t=1}^\tau \xi_t\left\langle e_t, x_t-x\right\rangle \leq \frac{\epsilon W_\tau}{r-\epsilon} .
$$
This relation holds true for every $x \in Q_1$, and \eqref{eq:resid2} follows.

\section{Proof of Theorem \ref{thm:D_lower} and Corollary \ref{cor:convergence}}\label{proof:D_lower}

The proof of Theorem \ref{thm:D_lower} is divided into three parts. First,
we ``lift'' the original space $\R^n$, treating it as a hyperplane $E=\{(x,s)\in\R^{n+1} \mid s=1\}$, and introduce a set $Q_{\tau+1}^{+}$. In the second part, we describe the polar of this set. Both $Q_{\tau+1}^{+}$ and its polar play an important role in the third part of the proof, where we provide a lower bound on the optimal value.


\subsection{``Lifting'' the Original Space}

Let us treat the original space $\R^n$ as a hyperplane in $\R^{n+1}$, that is, $E=\{(x,s)\in\R^{n+1} \mid s=1\}$. Define the set $Q_{\tau+1}^+\subset \R^{n+1}$ as a convex hull of the origin $0\in\R^{n+1}$ and $Q_{\tau+1}$ (treated as a subset of a hyperplane $E\subset \R^{n+1}$).
Let $\bar{A}:=\left[\begin{array}{cc}A & -b \end{array}\right]$ represent the constraints that form $Q_{\tau+1}$. We will now show that
\begin{equation}\label{eq:conv_hull}
    Q_{\tau+1}^+ = \left\{z \in \R^{n+1}: \bar{A}z \leq \nul \right\} \cap \left\{z \in \R^{n+1}: \bigl\la \last, z \bigr\ra \leq 1 \right\}.
\end{equation}

First note that since $Q_{\tau+1}$ is a bounded polytope, the system of inequalities $Ay \leq 0$ only has a trivial solution. Indeed, if we had a nonzero solution $y$, then for any $x\in Q_{\tau+1}$, the ray $x+\alpha y, \alpha\geq 0$ would belong to $Q_{\tau+1}$: $A(x+\alpha y) = Ax + \alpha Ay \leq b$, which contradicts the boundedness of $Q_{\tau+1}$.

Let $M$ be the right-hand side of \eqref{eq:conv_hull}.
Let us show that if $\left[\begin{smallmatrix} x \\ s \end{smallmatrix}\right] \in M$, then $s\geq 0$. Indeed, if $s<0$, then
\begin{equation*}
    Ax \leq sb \Rightarrow A\frac{x}{|s|} \leq -b.
\end{equation*}
Since $Q_{\tau+1}$ has a nonempty interior, there exists a point $y\in \R^n$ with $Ay < b$, therefore, $A(y+\frac{x}{|s|}) < b-b = 0$, which is impossible.

It is evident that $M$ is convex and contains both $0$ and $Q_{\tau+1}$ (treated as a subset of a hyperplane $E\subset \R^{n+1}$). What is left to prove is that any point $\left[\begin{smallmatrix} x \\ s \end{smallmatrix}\right] \in M$ is a convex combination of a point in $Q_{\tau+1}$ and $0$. If $s=0$, then $Ax\leq 0$ and we conclude that $x=0$. In the opposite case, we have $s\in (0,1]$. The vector $y:=s^{-1}x$ satisfies $Ay =s^{-1}Ax \leq b$ since $Ax \leq sb$, i.e., $\left[\begin{smallmatrix} y \\ 1 \end{smallmatrix}\right] \in Q_{\tau+1}$. Thus,
\begin{equation*}
    \left[\begin{smallmatrix} x \\ s \end{smallmatrix}\right] = s \left[\begin{smallmatrix} y \\ 1 \end{smallmatrix}\right] + (1-s) 0,
\end{equation*}
which concludes the proof of \eqref{eq:conv_hull}.

\subsection{Polar of a Set}

The polar of a set $P\subseteq \R^{n+1}$ is defined as
$$
    \pol P:=\left\{z \in \R^{n+1} \mid \la z, p\ra \leq 1\; \forall p \in P\right\}.
$$
We will now show that $\pol Q_{\tau+1}^+$ has the form
\begin{equation}\label{eq:polar}
    \pol Q_{\tau+1}^+ = \left\{ \left[\begin{smallmatrix} x \\ g \end{smallmatrix}\right] \in \R^{n+1} \mid x = A^\top \lambda,\, g=-b^\top \lambda + s,\, \lambda \in \R_+^{m},\, s \in [0,1] \right\}.
\end{equation}
To do so, we will use the following 
\begin{lemma}[\cite{nemirovski2010accuracy}, Lemma 6.3]
\label{Lm:nemirovski2010accuracy}
    Let $P, Q$ be two closed convex sets in $\R^{n+1}$ containing the origin and such that $P$ is a cone, and let $\int P \cap \int Q \neq \emptyset$. Then $\pol(P \cap Q)=\pol Q+P_*$, where
    $$
        P_*=\left\{z \in \R^{n+1}: \la z, u\ra \leq 0\; \forall u \in P\right\}.
    $$
\end{lemma}
Observe that $P:=\left\{z \in \R^{n+1}: \bar{A} z \leq \nul \right\}$ and $Q:=\left\{z \in \R^{n+1}: \bigl\la \last, z \bigr\ra \leq 1 \right\}$ are closed convex sets containing the origin, $P$ is a cone, and $\int P \cap \int Q \neq \emptyset$. Thus, the lemma applies.
Polar of a polyhedral cone $P$ is a finitely generated cone (see, for example, Lemma 1.12 (4) in \cite{paffenholz2010polyhedral}), i.e.,
\begin{equation}\label{eq:fin_gen_con}
    P_*=\left\{ z = \bar{A}^\top \alpha: \alpha \in \R_+^m \right\}.
\end{equation}

Let us show that 
\begin{equation}\label{eq:pol_halfspace}
    \pol Q = \left\{ s \last: s \in [0,1] \right\}.
\end{equation}
Denote the right hand side of \eqref{eq:pol_halfspace} by $\tilde{Q}$.
Let $y \in \tilde{Q}$, i.e., $y = s \last,\, s \in [0,1]$, then for any $z\in Q$, we have $\la y, z\ra = s \la \last, z \ra \leq s \leq 1 \Rightarrow y \in \pol Q$. Now let $y = \left[\begin{smallmatrix} x \\ s \end{smallmatrix}\right] \notin \tilde{Q}$,
i.e., $s \notin [0,1]$ or $x \neq \nul$. 
\begin{itemize}
    \item If $s<0$, then for $z = \left[\begin{smallmatrix} \nul \\ s^{-1}-1 \end{smallmatrix}\right]$ it holds $z \in Q$ and $\la y, z\ra = 1-s>1 \Rightarrow y \notin \pol Q$.
    \item If $s>1$, then for $z = \last$ it holds $z \in Q$ and $\la y, z\ra = s>1 \Rightarrow y \notin \pol Q$.
    \item If $x \neq \nul$, then for $z = \left[\begin{smallmatrix} 2x/\|x\|_2^2\\ 0 \end{smallmatrix}\right]$ it holds $z \in Q$ and $\la y, z\ra = 2 > 1 \Rightarrow y \notin \pol Q$.
\end{itemize}
Thus, $\pol Q = \tilde{Q}$ and the formula \eqref{eq:polar} follows from \eqref{eq:fin_gen_con}, \eqref{eq:pol_halfspace} and Lemma \ref{Lm:nemirovski2010accuracy}.

\subsection{Lower Bound on an Optimal Value}

Consider the ellipsoid $\mathcal{E}$ of maximal volume contained in $Q_{\tau+1}$. It is called John ellipsoid and it has a property that
\begin{equation}\label{eq:john}
    \mathcal{E} \subset Q_{\tau+1} \subset \hat{\mathcal{E}} := \{nx\mid x \in \mathcal{E}\},
\end{equation}
see, e.g., \cite{boyd2004convex}, Chapter 8.4.
As it was shown in \cite{nemirovski2010accuracy} (subsection 4.3), for an ellispoid $\hat{\mathcal{E}}$ there exists a vector $h'\in \R^n$ with $\|h'\|_2 \geq \frac{1}{2 \rho(\hat{\mathcal{E}})}$ such that
$$
\max_{x \in \hat{\mathcal{E}}}\langle h', x\rangle-\min_{x \in \hat{\mathcal{E}}}\langle h', x\rangle \leq 1.
$$
As it follows from \eqref{eq:john}, $\|h'\|_2 \geq \frac{1}{2n \rho(\mathcal{E})} \geq \frac{1}{2n \rho(Q_{\tau+1})}$ and
\begin{equation}\label{eq:narrowest_stripe}
   \max_{x \in Q_{\tau+1}}\langle h', x\rangle-\min_{x \in Q_{\tau+1}}\langle h', x\rangle \leq 1.
\end{equation}
The last formula implies
that both vectors
$$
h^{+}=\left[\begin{smallmatrix} h' \\ -\la h', x_{\tau+1} \ra \end{smallmatrix}\right],\;
h^{-}=-h^{+}
$$
belong to $\pol\left(Q_{\tau+1}^{+}\right)$ since for any $z\in Q_{\tau+1}^{+}$ it holds
$$
    z=\left[\begin{smallmatrix} sx \\ s \end{smallmatrix}\right]\quad \text{for some } x\in Q_{\tau+1},\, s\in[0,1],
$$
therefore,
\begin{align*}
    \la h^+, z\ra = \la h', sx\ra - \la h', x_{\tau+1}\ra \cdot s \leq \max_{x \in Q_{\tau+1}}\langle h', x\rangle-\min_{x \in Q_{\tau+1}}\langle h', x\rangle \stackrel{\eqref{eq:narrowest_stripe}}{\leq} 1,
\end{align*}
and similarly for $h^-$.
According to \eqref{eq:polar}, there exist $\mu,\eta \in \R^m_+$ and $u,v,\in [0,1]$ such that
\begin{equation}\label{eq:h_decomp}
    h^{+} = \left[\begin{smallmatrix} A^\top \mu \\ -b^\top \mu+u \end{smallmatrix}\right],\; h^{-} = \left[\begin{smallmatrix} A^\top \eta \\ -b^\top \eta+v \end{smallmatrix}\right].
\end{equation}
Observe that
\begin{equation*}
    0 = h^{+} + h^{-} = \left[\begin{smallmatrix} A^\top(\mu+\eta) \\ -b^\top(\mu+\eta)+u+v \end{smallmatrix}\right],\; u+v\in [0,2],
\end{equation*}
i.e., $\lambda := \mu+\eta$ is a feasible point for the LP problem \eqref{eq:lp}.

Let $\bar{x}$ be the center of Euclidean ball $B$ of radius $r$ which is contained in $X$. Consider first the case when $\left\langle h', \bar{x}-x_{\tau+1}\right\rangle \geq 0$. Multiplying $h^+$ by $x^{+}=[\bar{x}+r e ; 1]$ with $e \in \R^n,\|e\|_2 \leq 1$, we get
\begin{equation}\label{eq:h+x+0}
    \la h^{+}, x^{+}\ra = \la h', r e\ra + \la h', \bar{x}-x_{\tau+1}\ra \geq \langle h', r e\rangle.
\end{equation}
At the same time,
\begin{equation}\label{eq:h+x+}
\begin{aligned}
    \la h^{+}, x^{+}\ra & \stackrel{\eqref{eq:h_decomp}}{=} \Bigl\la \sum_{i=1}^m \mu_i a_i, \bar{x}+r e \Bigr\ra - \sum_{i=1}^m \mu_i b_i + u \\
    &\leq \sum_{i \in \Pt \cup \Nt} \mu_i \la a_i, \bar{x}+r e-x_{t(i)}\ra + \sum_{i \in \It} \mu_i \left(a_i^\top(\bar{x}+re)-b_i\right) + u,
\end{aligned}
\end{equation}
where we used $a_i^\top x_{t(i)}\leq b_i,\, i \in \Pt \cup \Nt$ (see the end of subsection \ref{subsec:polytope_algo}). Further, since $u\leq 1$ and $\bar{x}+re \in Q_1 \Rightarrow a_i^\top(\bar{x}+re)\leq b_i, i\in \It$, it holds
\begin{equation}
\begin{aligned}\label{eq:h+x+2}
    \la h^{+}, x^{+}\ra &\stackrel{\eqref{eq:h+x+}}{\leq} \sum_{i \in \Pt \cup \Nt} \mu_i\la a_i, \bar{x}+r e-x_{t(i)}\ra + 1, \\
    &\leq \sum_{i \in \Pt} \mu_i \la a_i, \bar{x}+r e-x_{t(i)} \ra + 1,
\end{aligned}
\end{equation}
where the last inequality is due to the fact that $a_i$ separates $x_{t(i)}$ and $X$ for $i \in \Nt$ which means $a_i^\top x \leq a_i^\top x_{t(i)}$ for all $x\in X$.
Note that $x_{t(i)}\in \int X\subseteq Q_1$ for all $i\in \Pt$. Thus, combining \eqref{eq:h+x+0} and \eqref{eq:h+x+2}, we obtain
$$
    \langle h', r e\rangle \leq 1+\sum_{i \in \Pt} \mu_i \left\|a_i\right\|_2 D\left(Q_1\right) \leq 1 + D_\tau D\left(Q_1\right).
$$
The resulting inequality holds true for all unit vectors $e$; maximizing the left hand side over these $e$, we get $D_\tau \geq \frac{r\|h'\|_2-1}{D\left(Q_1\right)}$. Recalling that $\|h'\|_2 \geq \frac{1}{2 n \rho\left(Q_{\tau+1}\right)}$, we arrive at \eqref{eq:D_lower}. We have established it in the case of $\left\langle h', \bar{x}-x_{\tau+1}\right\rangle \geq 0$; in the opposite case we can use the same reasoning with $h^-$ in the role of $h^+$.

\subsection{Proof of Corollary \ref{cor:convergence}}\label{proof:convergence}

The certificate is well-defined since
\begin{equation}\label{eq:small_rho}
    \rho\left(Q_{\tau+1}\right) \leq \frac{(1-\alpha) r^2}{16 n D\left(Q_1\right)} < \frac{r}{2n}
\end{equation}
implies 
\begin{equation*}
    D_\tau(\lambda) \geq (1-\alpha) D_\tau(\lambda^*) \stackrel{\text{Thm} \ref{thm:D_lower}}{\geq} (1-\alpha) D^{-1}\left(Q_1\right)\left(\frac{r}{2 n \rho\left(Q_{\tau+1}\right)}-1\right) \stackrel{\eqref{eq:small_rho}}{>} 0,
\end{equation*}
therefore, $\lambda_i>0$ for some $i\in \Pt\Rightarrow d_\tau>0$.

Let us now show that $\epsilon_\tau := \frac{2}{D_\tau}<r$.
\begin{align*}
    \rho\left(Q_{\tau+1}\right) &\leq \frac{(1-\alpha) r^2}{16 n D\left(Q_1\right)} = \frac{r}{2n}\left( \frac{8D\left(Q_1\right)}{(1-\alpha) r} \right)^{-1} < \frac{r}{2n}\left( \frac{2D\left(Q_1\right)}{(1-\alpha) r}+1 \right)^{-1} \\
    & \Rightarrow (1-\alpha) D^{-1}\left(Q_1\right)\left(\frac{r}{2 n \rho\left(Q_{\tau+1}\right)}-1\right) > \frac{2}{r} \\
    & \Rightarrow D_\tau(\lambda) \geq (1-\alpha) D_\tau(\lambda^*) \stackrel{\text{Thm} \ref{thm:D_lower}}{\geq} (1-\alpha) D^{-1}\left(Q_1\right)\left(\frac{r}{2 n \rho\left(Q_{\tau+1}\right)}-1\right) > \frac{2}{r}.
\end{align*}
Applying Lemma \ref{lem:resid}, we get
\begin{align*}
    \resid\resar &\leq \frac{\epsilon_\tau}{r-\epsilon_\tau} W_\tau = \left( \frac{r D_\tau}{2}-1 \right)^{-1} W_\tau \\
    &\stackrel{\text{Thm} \ref{thm:D_lower}}{\leq} \left( \frac{(1-\alpha) r}{2 D\left(Q_1\right)}\left(\frac{r}{2 n \rho\left(Q_{\tau+1}\right)}-1\right) - 1 \right)^{-1} W_\tau \\
    &\leq \left( \frac{(1-\alpha) r^2}{8 n \rho\left(Q_{\tau+1}\right) D\left(Q_1\right)} - 1 \right)^{-1} W_\tau,
\end{align*}
where the last inequality follows from
$$
    \rho\left(Q_{\tau+1}\right) \leq \frac{r}{4n} \Rightarrow \frac{r}{2 n \rho\left(Q_{\tau+1}\right)}-1 \geq \frac{r}{4 n \rho\left(Q_{\tau+1}\right)}.
$$
Finally, using
$$
    \rho\left(Q_{\tau+1}\right) \leq \frac{(1-\alpha) r^2}{16 n D\left(Q_1\right)} \Rightarrow \frac{(1-\alpha) r^2}{8 n \rho\left(Q_{\tau+1}\right) D\left(Q_1\right)} - 1 \geq \frac{(1-\alpha) r^2}{16 n \rho\left(Q_{\tau+1}\right)},
$$
we arrive at the first result. The second result now follows directly from the definition of $\delta$-subgradient \eqref{eq:delta_subgrad}.


\bibliographystyle{spmpsci}
\bibliography{references}









\end{document}